\documentclass[11pt,twoside]{article}
\usepackage{amsmath,amsbsy,amsfonts}
\usepackage[notref,notcite]{}
\usepackage{hyperref}
\usepackage{enumerate}
\usepackage{enumitem}
\usepackage{pdfsync}
\usepackage{color}
\usepackage{amssymb}

\usepackage[latin1]{inputenc}
\newtheorem{lem}{Lemma}[section]
\newtheorem{prop}{Proposition}[section]
\newtheorem{thm}{Theorem}[section]

\newtheorem{rmk}{Remark}[section]
\newtheorem{claim}{Claim}[section]
\newtheorem{defi}{Definition}[section]

\let\Section=\section

\def\section{\setcounter{equation}{0}\Section}

\def\nd{\noindent}

\def\hsf{\hspace*{\fill}}
\def\proof{{\rm \bf Proof}}

\newcommand{\w}{W_{0}^{1}L_{\Phi}(\Omega)}

\def\r{\mathrm{R}}

\begin{document}

\title{ On Strongly Nonlinear  Eigenvalue Problems in the Framework of  Nonreflexive Orlicz-Sobolev Spaces}
\vskip.2cm

\author{Edcarlos D. Silva~~~~Jose V. A. Goncalves\\
\\
Kaye  O.  Silva}
\vskip.7cm

\date{}

\pretolerance10000

\maketitle

\begin{abstract}
\noindent {}
\vskip.5cm
It is established  existence and multiplicity of solutions for strongly nonlinear   problems driven by the $\Phi$-Laplacian operator on bounded domains.  Our main results are stated without the so called $\Delta_{2}$ condition at infinity which means that the underlying Orlicz-Sobolev spaces are not reflexive.\\
\\

{\small

}
\end{abstract}

\section{Introduction}

\nd  In this work, we study the nonlinear eigenvalue problem
~\begin{equation}\label{3.1}
\left\{\begin{array}{rllr}
-{\rm div} \big( \phi(|\nabla u|)\nabla u\ \big ) &=& \lambda f(u)~~\mbox{in}~~\Omega,\\

u&=&0~~ \mbox{on}~~  \partial {\Omega,}
 \end{array}
\right.
\end{equation}
\nd where $\Omega\subset\r^N$ is a bounded domain with smooth boundary $\partial \Omega$,~$\lambda > 0$ is a parameter and $\phi : (0, \infty) \to (0, \infty)$ is a  $C^1$-function satisfying
\begin{itemize}
  \item[($\phi_1$)]  \     \ $\mbox{(i)} \ \ t\phi(t)\to 0 \ \mbox{as} \  t\to 0, \\
\\
 ~~ \mbox{(ii)} \  t\phi(t)\to\infty \ \mbox{as} \  t\to\infty$,
  \item[($\phi_2$)]  \ $t\phi(t) \ \mbox{\it is strictly increasing in}~ (0, \infty)$.
 \end{itemize}
\nd Throughout this work,  $f : [0, \infty) \to [0, \infty)$ is a continuous function satisfying
\begin{itemize}
  \item[($f_1$)] \  $f(0)\geq 0$,
  \item[($f_2$)] \  {\it there exist positive  numbers} $a_k, b_k,\ k = 1, \cdots, m-1$ {\it such that}
$$
0<a_1<b_1<a_2<b_2<...<b_{m-1}<a_m, \\
$$
$$
f(s)\leq 0 \ \mbox{ if } \ s \in (a_k,b_k), \\
$$
$$
f(s)\geq 0 \ \mbox{ if } \ s \in (b_k,a_{k+1}),
$$
\item[($f_3$)] \  $\displaystyle\int_{a_{k}}^{a_{k+1}}f(s)ds>0$,  $k = 1,  \cdots,  m-1$.
\nd
\end{itemize}

\begin{rmk}
We  extend   $t \mapsto t \phi(t)$ to the whole of $\r$ as an odd function.
\end{rmk}

\nd We shall consider the N-function
$$
\Phi(t) = \int_0^t s \phi(s) ds,~~ t \in \r,
$$
\nd and we shall use the notation
$$
 \Delta_{\Phi}  u  :={\rm div} \left( \phi(|\nabla u|)\nabla u\ \right)
$$
\nd for the  $\Phi$-Laplacian operator. Due to the use of this more general operator we shall work in the framework of Orlicz-Sobolev spaces such as $\w$ which under our assumptions on $\phi$ is not reflexive.
\vskip.1cm

\nd The main novelty in this work is to ensure existence and multiplicity of solutions for problem \eqref{3.1} without the so called $\Delta_{2}$ condition which  amounts  that $\w$ is not reflexive.
\vskip.1cm

\nd Examples of functions $\Phi$ covered by the main theorem  in the present work:
\begin{equation}\label{ex 3}
\Phi(t)=e^t -t+1,
\end{equation}
\begin{equation}\label{ex 1}
\Phi(t)=(1+t^2)^{\gamma} - 1~~ \mbox{where}~~\gamma>\frac{1}{2},
\end{equation}
\begin{equation}\label{ex 2}
\Phi(t)=t^p \log(1+t)~~ \mbox{where}~~p\geq 1.
\end{equation}

\nd Our main result stated below extends Theorem 1.1 by Loc \& Schmitt in \cite{loc-schmitt} to the more general operator $\Delta_{\Phi}$ in the case that $W_{0}^{1}L_{\Phi}(\Omega)$ is not reflexive.
\vskip.1cm

\begin{thm}\label{1.1}
Assume  $(\phi_1)-(\phi_2)$.  Then
\vskip.1cm
\nd {\rm (i)}~ if~ $(f_1)-(f_3)$ hold,  there is $\overline{\lambda}>0$  such  that for each $\lambda>\overline{\lambda}$, $(\ref{1.1})$ admits at least $m-1$ non negative weak solutions  $u_1,...,u_{m-1} \in \w \cap L^{\infty}(\Omega)$ satisfying
$$
a_1 <  \|u_{1}\|_\infty \leq a_2 < \|u_{2}\|_\infty \leq \cdots \leq a_{m-1}<\|u_{m-1}\|_\infty \leq a_{m},
$$
\nd {\rm (ii)}~conversely, if $u \in \w \cap L^{\infty}(\Omega)$ is a nonnegative weak solution of problem $(\ref{3.1})$ such that $a_k < \|u\|_\infty \leq a_{k+1}$ and $(f_1)-(f_2)$ holds then $(f_3)$  also holds true.
\end{thm}

\begin{rmk} If $f$ is extended to the whole of $\r$ as an odd function then with minor modifications on the arguments of the present work, problem \ref{3.1} admits at least $2(m-1)$ weak solutions $u_{1}, \ldots, u_{m-1}$ and $v_{1}, \ldots, v_{m-1}$ such that $u_{i} > 0, v_{i} < 0$ in $\Omega,i = 1, 2, \ldots, m -1$ and
$$
a_1 <  \|u_{1}\|_\infty \leq a_2 < \|u_{2}\|_\infty \leq \cdots \leq a_{m-1}<\|u_{m-1}\|_\infty \leq a_{m},
$$
$$
a_1 <  \|v_{1}\|_\infty \leq a_2 < \|v_{2}\|_\infty \leq \cdots \leq a_{m-1}<\|v_{m-1}\|_\infty \leq a_{m}.
$$
\end{rmk}


\nd We recall that Hess in \cite{Hess} employed variational and topological methods and arguments with lower and upper solutions to prove  a result on existence of multiple positive solutions for the problem
$$
- \Delta  u = \lambda f(u)~\mbox{in}~\Omega,~~u = 0~\mbox{on}~\partial \Omega,
$$
\nd which is problem \eqref{3.1} with  $\phi(t)  \equiv 1$. As noticed by  Hess, the results in \cite{Hess} were motivated by Brown \& Budin \cite{BB-1, BB-2} which in turn were motivated by  the literature on nonlinear heat generation.
\vskip.2cm

\nd In  \cite{loc-schmitt}, Loc \& Schmitt extended the result by Hess to the $p$-Laplacian operator by taking $\phi(t) = t^{p-2} \ \mbox{with} \  1 < p < \infty$ in (1.1). Actually, in \cite{loc-schmitt}  the authors showed  that $(f_1)-(f_3)$ are sufficient conditions for the existence of $m-1$ positive solutions of
$$
- \Delta_p  u = \lambda f(u)~\mbox{in}~\Omega,~~u = 0~\mbox{on}~\partial \Omega,
$$
\nd for $\lambda$ large,  while   if $(f_1)-(f_2)$ hold, and $u$ is a positive solution of the problem above
with $a_k < \|u\|_\infty \leq a_{k+1}$ then $(f_3)$ holds.
\vskip.1cm

\nd Regarding the rich literature on \cite{Hess,loc-schmitt}  we further refer the reader to  \cite{chipot, correa-3, Oliveira,dancer-schmitt}, where several techniques were employed.
\vskip.1cm

\nd In particular, in \cite{Oliveira}, the authors  proved a version of Theorem \ref{1.1}  in the case that both  $\Phi$ and its conjugate function $\tilde{\Phi}$ satisfy  the $\Delta_2$ condition
, or equivalently,
$$
\mbox{there exist}~ \ell, m~\mbox{with}~ 1< \ell < m < \infty~\mbox{such that}
$$
$$
\ell\leq\frac{t\Phi'(t)}{\Phi(t)}\le m,~ t>0,
$$
\nd see e.g. Section  \eqref{orlicz} for clarification on  notation and terms  above.
\vskip.2cm

\nd Examples of functions $\Phi$ for which the $\Delta_2$ condition does not hold are (\ref{ex 3}) and (\ref{ex 2}) with $p=1$. In these two cases $\Phi$ grows too slow or too fast, respectively. In the first case  we have that $\ell=1$ while in the second case $m=\infty$.
\vskip.1cm

\nd In the  present work, we do not require  the $\Delta_2$ condition. Due to the lack of the $\Delta_{2}$ condition, we have to overcome many difficulties which do not appear in the $\Delta_2$ case. Indeed, without that condition, the Orlicz spaces are not reflexive and the energy functional $J$ associated to problem (\ref{3.1}) is not $C^{1}$, in fact, it is not even well defined in the whole Orlicz-Sobolev space. This difficulty is overcome by working in some appropriate subspaces in the Orlicz-Sobolev space. On this subject we refer the reader to  Gossez \cite{Gz1} and  Garc\'ia-Huidobro et al \cite{huidobro}.
For further results without $\Delta_{2}$, we refer the reader to Le Vy Khoi \cite{Khoi}, Loc
\verb"&" Schmitt \cite{loc-schmitt}, V. Mustonen \& M. Tienari \cite{ Tienari} and
 M. Tienari \cite{Tienaritese}, Cl\'ement et. al. \cite{clem},
\vskip.1cm

\nd Other classes of functions $\phi$ which satisfy $(\phi_1)-(\phi_2)$ are:
\vskip.1cm
\begin{description}
\item{$\rm {(i)}$} $\phi(t)=t^{p-2}+t^{q-2}$ with $1<p<q<N$. In this case  with $\ell=p$ and $m=q$, the corresponding  operator is the $(p,q)$-Laplacian    and  \eqref{3.1}
   becomes
\begin{eqnarray}
\left\{\
\begin{array}{c}
\displaystyle-\Delta_p u-\Delta_q u= \lambda f(u)~\mbox{in}~\Omega,\\
u>0~\mbox{in}~\Omega,~u=0~\mbox{on}~\partial \Omega,
\end{array}
\right.\nonumber
\end{eqnarray}

\item{$\rm {(ii)}$} $\phi(t)=\displaystyle \sum_{i=1}^N t^{p_i-2}$ where $1<p_1<p_2<...<p_N$, $\frac{1}{\overline{p}}=\displaystyle \frac{1}{N}\sum_{i=1}^N\frac{1}{p_i}$ with $\overline p<N$. In this case the corresponding problem
\begin{eqnarray}
\left\{\
\begin{array}{c}
\displaystyle-\sum_{i=1}^N\Delta_{p_i} u= \lambda f(u)~\mbox{in}~\Omega,\\
u>0~\mbox{in}~\Omega,~u=0~\mbox{on}~\partial \Omega
\end{array}
\right.\nonumber
\end{eqnarray}

\nd is known in the literature as an anisotropic elliptic problem,

\item{$\rm {(iii)}$} $\phi(t)=a(t^p)t^{p-2}$ where $2\leq p<N$ and $a: (0,\infty) \to (0, \infty)$ is a suitable  $C^1(\mathbb{R}^+)$-function. In this case the corresponding problem reads as
\begin{eqnarray}\label{ex-geral}
\left\{\
\begin{array}{c}
\displaystyle-\mbox{div}(a(|\nabla u|^p)|u|^{p-2}\nabla u) = \lambda f(u)~\mbox{in}~\Omega,\\
u>0~\mbox{in}~\Omega,~u=0~\mbox{on}~\partial \Omega.
\end{array}
\right.
\end{eqnarray}
\end{description}

\nd This work is organized as follows: in Section 2 we present some tools and references on Orlicz-Sobolev spaces. Section 3 is devoted to some auxiliary problems and the variational setting. Section 4 is devoted to some technical lemmata. In Section 5 we give the proof of Theorem \ref{1.1}. In the Appendix we give, for completeness, some technical results used in this work.

\section{Basics on Orlicz-Sobolev Spaces}\label{orlicz}

\nd The main references for this Section are Gossez \cite{Gz1, Gz3}, Adams \cite{A}, Kufner \cite{Kuf} and Tienari \cite{Tienaritese}. Further references will be given timely.
\vskip.2cm

\nd Let $\Phi: \r \rightarrow \r$ be an $N$-Function. We say that $\Phi$ satisfies the $\Delta_{2}$ condition if there exist $t_{0} > 0$ and $K > 0$  such  that
\begin{equation}
\Phi(2 t) \leq K \Phi(t),\ |t| \geq t_{0}.
\end{equation}
\nd It is  well known  that this condition is equivalent to
$$
\frac{t\Phi'(t)}{\Phi(t)}\le m,\  |t|\ge t_0~ \mbox{for some}~ m \in (1,\infty).
$$
\nd The $\Delta_{2}$ condition is crucial to ensure that Orlicz and Orlicz-Sobolev spaces are reflexive Banach
spaces.
\vskip.1cm

\nd   The Orlicz class associated with $\Phi$  is
$$
\mathcal{L}(\Omega)=\left\{u:\Omega\to\mathbb{R}~|~  u\ \mbox{measurable and}\ \int_\Omega \Phi(u) dx<\infty\right\}.
$$
\nd  The Orlicz space $L_\Phi(\Omega)$ is  the linear hull of $\mathcal{L}_\Phi(\Omega)$, that is,
$$
L_\Phi(\Omega) = \displaystyle  \bigcap_{{\mathcal{L}}(\Omega) \subset V }  \{V~|~  V \mbox{is a vector space}~   \}.
$$
\nd   The  norm of a function $u\in L_\Phi(\Omega)$, (Luxemburg norm),  is defined by
$$
\|u\|_\Phi=\inf\left\{\lambda > 0 ~|~\ \int_\Omega\Phi\left(\frac{u}{\lambda}\right)\le 1\right\}
$$
\nd and $L_\Phi(\Omega)$ endowed with the norm  $\|\cdot\|_\Phi$ is a Banach space. On the other hand, the space $E_\Phi(\Omega)$ is  defined by
$$
E_\Phi(\Omega) =   \displaystyle  \mbox{closure of } L^{\infty}(\Omega)~\mbox{in}~  L_\Phi(\Omega)~ \mbox{with respect to}~\Vert \cdot \Vert_{\Phi}.
$$
\nd We recall that
$$
L_\Phi(\Omega)= \displaystyle \Big \{u~|~ \int_\Omega \Phi(\lambda u)<\infty\ \mbox{for some}\ \lambda>0 \Big \}.
$$
\nd and
$$
E_\Phi(\Omega)= \Big \{u~|~ \int_\Omega\Phi(\lambda u)<\infty~ \mbox{for each}\ \lambda>0 \Big \}.
$$
\nd It is well known that $L_{\Psi}(\Omega)\hookrightarrow L^1(\Omega)$ (cf. Adams \cite{A}).
\vskip.1cm

\nd The Orlicz-Sobolev space is defined by
$$
 W^1L_\Phi(\Omega) = \displaystyle \Big\{u \in L_\Phi(\Omega)~|~ \frac{\partial u}{\partial x_i} \in L_\Phi(\Omega),~ i=1,...,N \Big\},
$$
\nd and similarly
$$
 W^1E_\Phi(\Omega)= \displaystyle \Big\{u \in E_\Phi(\Omega)~|~ \frac{\partial u}{\partial x_i} \in E_\Phi(\Omega),~ i=1,...,N \Big\}
$$
\nd  It is known  that $W^1L_\Phi(\Omega)$ endowed with the norm
$$
\|u\|_{1,\Phi}=\|u\|_\Phi+\sum_{i=1}^N \left\|\frac{\partial u}{\partial x_i}\right\|_{\Phi}
$$
\nd is a Banach space
\vskip.1cm

\nd We point out that  by the Poincar\'e Inequality  (see \cite[Lemma 5.7 p 202]{Gz1}),
$$
\|u\| :=\|\nabla u\|_\Phi,~u \in W^1L_\Phi(\Omega)
$$
\nd defines a norm in $W_0^1L_\Phi(\Omega)$ equivalent to  $\|\cdot \|_{1,\Phi}$.
\vskip.2cm

\nd The conjugate function of $\Phi$ is defined by
$$
\widetilde{\Phi}(t)=\sup\{ts-\Phi(s)~|~ s\in\mathbb{R}\}.
$$
\nd  It is known  that $\widetilde{\Phi}$ is also an   N-function and actually
$$
\widetilde{\Phi}(t)=\int_0^{|t|} t\widetilde{\phi}(t)dt~ \mbox{for a suitable  function}~  \widetilde{\phi}.
$$
\nd where  $\widetilde{\phi}$ satisfies the same basic conditions as $\phi$.  The Young inequality holds,
\begin{equation}\label{young}
ts\le \Phi(t)+\widetilde{\Phi}(s),\ t,s\in\mathbb{R},
\end{equation}
\nd and the equality is true if and only if $t=s\widetilde{\phi}(s)$ or $s=t\phi(t)$.
\vskip.2cm

\nd  It is well known  that $L_\Phi(\Omega)$ is the dual space of $E_{\widetilde{\Phi}}(\Omega)$, that is

$$
L_\Phi(\Omega) =  E_{\widetilde{\Phi}}(\Omega)^{\prime}.
$$
\nd  The H\"older inequality is true,  that is
 $$
\int_\Omega |uv|\le 2\|u\|_\Phi\| v\|_{\widetilde{\Phi}},~ u\in L_\Phi(\Omega),\ v\in L_{\widetilde{\Phi}}(\Omega).
$$
\nd  Yet following  Gossez \cite{Gz1, Gz3} we recall that the spaces $W^1L_\Phi(\Omega)$ and $W^1E_\Phi(\Omega)$ are  identified with  subspaces of the product spaces  $\prod L_\Phi(\Omega)$,~ $\prod E_\Phi(\Omega)$, respectively.
\vskip.2cm

\nd  We define
$$
W_0^1L_\Phi(\Omega) = \displaystyle {{\overline{{C_{0}^{\infty}}(\Omega)}}}^{\sigma\left(\prod L_\Phi(\Omega),\prod E_{\widetilde{\Phi}}(\Omega)\right)},
$$
$$
W_0^1E_\Phi(\Omega) = \displaystyle  {\overline{C_{0}^{\infty}(\Omega)}}^{(E_{\Phi}, \|u\|_{1,\Phi})}.
$$

\section{ Auxiliary Problems and Variational Setting}

\nd Consider the family of problems associated to (\ref{3.1})
~\begin{equation}\label{3.1k}
\left\{\begin{array}{rllr}
-\Delta_{\Phi} u &=& \lambda f_k(u)~~\mbox{in}~~\Omega,\\

u&=&0~~ \mbox{on}~~  \partial {\Omega,}
 \end{array}
\right.
\end{equation}
\nd where $f_k : \r \to \r$ is a continuous function for each $k = 2, \cdots, m$ which is given by
$$
f_k(s) = \left\{ \begin{array}{rl}
 f(0) &\mbox{if} s\leq0 ,\\
  f(s) &\mbox{if}~ 0\leq s\leq a_k , \\
  0 &\mbox{if}~ s>a_k.
       \end{array} \right.
$$
Here we emphasize that $f_{k}$ is a bounded function which has at least $m$ bumps.
\nd The energy functional  associated to the problem $(\ref{3.1k})$ is given by
$$
I_k(\lambda,u)=\int_\Omega\Phi(|\nabla u|) dx-\lambda\int_\Omega F_k(u) dx,~ u \in W_0^1L_{\Phi}(\Omega),
$$
\nd where
$$
F_k(s)=\int_0^s f_k(t)dt, k = 1, 2, \dots, m -1.
$$
\nd  Now we observe that
$I_k(\lambda,\cdot):W_0^1L_\Phi(\Omega)\to\r\cup\{\infty\}$. In fact, the domain of $I_k$, that is the set for which $I_k$ is finite, is
$$
\{u\in W_0^1L_\Phi(\Omega)~|~  |\nabla u|\in \mathcal{L}(\Omega)\}.
$$
 Moreover, there are points in $W_0^1L_{\Phi}(\Omega)$ where $I_k$ is not differentiable.  However, we mention that $I_k$ is differentiable in $W_0^1E_{\Phi}(\Omega)$ , (cf. \cite[Lemma 3.4]{huidobro}). In this way, we say that $u\in W_0^1L_{\Phi}(\Omega)$ is a weak solution for problem $(\ref{3.1k})$ if
$$
\int_\Omega\phi(|\nabla u|)\nabla u \cdot \nabla v   dx = \lambda\int_\Omega f_k (u)v dx,\  v\in W_0^1L_\Phi(\Omega).
$$
\begin{rmk}
 When working with the $\Delta_2$ condition in both $\Phi$ and $\tilde{\Phi}$, the energy functional $I_k$ is $C^1$ (which is easy to prove), and so, every critical point of $I_k$ satisfies the Euler equation above. Without $\Delta_2$ condition, the lack of differentiability is a problem which we have to overcome in order to show that the minimum of $I_k$ satisfies the above Euler equation.
\end{rmk}

\section{Technical Lemmas}

\nd The result below is crucial in this work and it was proved originally by Hess \cite{Hess} and then extended by Loc \& Schmitt in \cite{loc-schmitt} to more general Sobolev spaces. In the present paper due to the non-reflexivity of the Orlicz-Sobolev space $\w$ it was necessary  to employ a version of Stampacchia's theorem for that space. See Proposition \ref{CA1}in the Appendix.
\vskip.1cm

\nd Consider the problem
\begin{equation}\label{2.1}
 \left\{ \begin{array}{c}
 -\Delta_\Phi u= g(u) ~\mbox{ in}~~ \Omega, \\
u   = 0~\mbox{on}~\partial \Omega .
       \end{array} \right.
\end{equation}
\nd The result is:

\begin{lem}\label{lema 2.1}
 Let $g:\r\to\r$ be a continuous function such  that $g(s)\geq 0$ for $s\in (-\infty,0)$ and assume that there is some $s_0\geq 0$  such that $g(s)\leq 0$ for $s\geq s_0$.  Let $u \in \w$ be a weak solution for the problem \eqref{2.1}.
\nd Then $0 \leq u \leq s_0~\mbox{a.e. in}~\Omega$.
\end{lem}

\nd \proof. Let $\Omega_0=\{x\in\Omega:\ u(x)<0\}$. Using Proposition \ref{CA1} in the Appendix, with the function $\min\{t,0\}$, we can take $u^-$ as a test function. As a consequence we know that
$$
\int _{\Omega_0}\phi(|\nabla u|)|\nabla u|^2 dx =\int _{\Omega_0}g(u)udx.
$$

\nd In particular, the last assertion says that
$$
\int _{\Omega_0}\phi(|\nabla u|)|\nabla u|^2 dx \le 0,
$$

\nd which implies that $\Omega_0$ has zero measure. Again, define $\Omega_{s_0}=\{x\in \Omega:\ u(x)>s_0\}$. Using Proposition \ref{CA1}, with the function $\max\{t-s_0,0\}$, we take $(u-s_0)^+$ as a test function. As a byproduct we get
$$
\int_{\Omega_{s_0}}\phi(|\nabla u|)|\nabla u|^2 dx =\int _{\Omega_{s_0}}g(u)(u-s_0) dx,
$$
\nd showing that
$$
\int_{\Omega_{s_0}}\phi(|\nabla u|)|\nabla u|^2 dx \le 0.
$$
\nd  Therefore the set $\Omega_{s_0}$ has zero measure. \hfill \fbox \hsf
\vskip.3cm

\nd  The result below is crucial. In the case that the $\Delta_2$ condition holds its proof is rather straightforward. In the setting of the present paper it is much more difficult. We shall detail it by  using arguments employed in \cite{huidobro} and \cite{Tienari}.
	
\begin{lem}\label{Qk not-empty}
Let $\lambda > 0$. Then there is $v_k \equiv v_k(\lambda) \in \w$  such that
$$
I_{k}(\lambda, v_k) = \min_{u \in \w } I_{k}(\lambda, u).
$$
\end{lem}

\nd   \proof . It is enough to show that $I_{k}(\lambda, \cdot)$ is both coercive and weak* sequentially lower semicontinuous, ($w^{\star}.$s.l.s.c.  for short).
\vskip.2cm

\nd In order to show  the coercivity property for $I_{k}(\lambda, )$  take $u \in \w$ such that $\|u\|\ge 1+\epsilon$ with $\epsilon>0$. Once $\Phi$ is convex and satisfies $\Phi(0)=0$ we have that
$$
\frac{1+\epsilon}{\|u\|}\int_\Omega \Phi(|\nabla u|)dx\ge \int_\Omega \Phi\left(\frac{(1+\epsilon)|\nabla u|}{\|u\|}\right) dx >1.
$$

\nd In the last inequality it was used the definition for Luxemburg norm. Here we point out that $f_{k}$ is a bounded function for any $k = 1, 2, \ldots, m-1$.  The coerciveness follows by a straightforward argument.
\vskip.2cm

\nd In order to show that
$I_{k}(\lambda, \cdot)$~ is $w^{\star}.$s.l.s.c.,  we first prove that
$$
 u \in \w \mapsto \int_\Omega \Phi(|\nabla u|) dx
$$
\nd is $w^{\star}.$s.l.s.c.. Indeed, it follows by Young's inequality ($\ref{young}$) that
$$\int_\Omega \Phi(|\nabla u|)dx =\sup\left\{\int_\Omega |\nabla u| w dx-\int_\Omega \tilde{\Phi}(w)dx~|~ w\in E_{\tilde{\Phi}}(\Omega)\right\}.$$

\nd Assume that $u_n \stackrel {w^{\star}} \rightharpoonup u$. Given $\epsilon >0$ it follows from the previous identity that there is $w\in E_{\tilde{\Phi}}(\Omega)$ such that
$$
\int_\Omega \Phi(|\nabla u_n|) dx \ge \int_\Omega |\nabla u_n|w dx-\int_\Omega \tilde{\Phi}(w) dx
$$
\nd and
$$\int_\Omega \Phi(|\nabla u|) dx \le \int_\Omega |\nabla u|w dx -\int_\Omega \tilde{\Phi}(w) dx +\epsilon.$$

\nd Due the fact that $|\nabla u|,|\nabla u_n|\ge 0$, we may assume without loss of generality that $w\ge 0$. As a consequence
\begin{eqnarray}\label{weaks}
 \int_\Omega \Phi(|\nabla u_n|)dx -\int_\Omega \Phi(|\nabla u|)dx &\geq& \int_\Omega |\nabla u_n|w dx-\int_\Omega |\nabla u|w dx-\epsilon \nonumber \\
&=& \int_\Omega |\nabla u_n w|dx-\int_\Omega |\nabla u w|dx-\epsilon.
\end{eqnarray}

\nd Since
$$
\int_\Omega \frac{\partial u_n}{\partial x_i}v dx\to \int_\Omega \frac{\partial u}{\partial x_i}v dx,~ v\in E_{\tilde{\Phi}}(\Omega),
$$

\nd and $E_{\widetilde{\Phi}}(\Omega) L^\infty(\Omega)=E_{\widetilde{\Phi}}(\Omega)$ (see \cite{Tienaritese}), we have that $\nabla u_n w\to \nabla u w$ for $\sigma(\prod L^1(\Omega),\prod L^\infty(\Omega))$. Hence, by the weak lower semicontinuity of norms, we conclude that
$$
\int_\Omega |\nabla u|w dx\le \liminf \int_\Omega |\nabla u_n|w dx.
$$

\nd This inequality and ($\ref{weaks}$) imply that
$$
\liminf \int_\Omega \Phi|\nabla u_n|) dx\ge \int_\Omega \Phi(|\nabla u|) dx-\epsilon.
$$

\nd Since $\epsilon$ was taken arbitrarily, we obtain that
$$
u \mapsto \int_\Omega \Phi(|\nabla u|) dx, u \in \w
$$
\nd  is a $w^{\star}$.s.l.s.c. function. Now, we shall prove that
$$
\int_\Omega F(u_n) dx \to \int_\Omega F(u) dx
$$
\vskip.3cm
 \nd First of all,  using the compact embedding $W_0^1L_{\Phi}(\Omega)\stackrel {\scriptsize cpt}   \hookrightarrow  L_{\Psi} (\Omega)$ for $\Psi << \Phi^\star$ and the continuous embedding $L_{\Psi}(\Omega)\hookrightarrow L^1(\Omega)$  we conclude that
$$
W_0^1L_{\Phi}(\Omega) \stackrel {\scriptsize cpt}   \hookrightarrow  L^1 (\Omega).
$$
\nd  Consequently we can take a function $h \in L^1(\Omega)$ such that
$$
 u_n \to u~ \mbox{and}~~  |u_n|\leq h~~ \mbox{a.e. in}~\Omega.
$$
\nd Using the dominated convergence theorem we have
$$
\int_\Omega F_k(u_n) dx \to \int_\Omega F_k(u) dx.
$$
\nd It follows from the previous arguments that
$$
I_{k}(\lambda, u) \leq \liminf I_{k}(\lambda, u_n).
$$
\nd As a consequence, there is a minimum $v_k \equiv v_k(\lambda)$ of $I_{k}(\lambda, \cdot)$. \hfill \fbox \hsf
\vskip.2cm
\nd Now we shall prove that $v_k$ is in fact a weak solution of problem ($\ref{3.1k}$). In order to do that, we will show that the functional $I_k$ is differentiable at its minimum point, found in Lemma \ref{Qk not-empty}. We first give some auxiliary results and definitions which can be found in \cite{Tienari}.
\vskip.2cm
\begin{defi} Let $$\operatorname{dom}(\phi(t)t)=\{u\in L_\Phi(\Omega):\ \phi(|u|)|u|\in L_{\tilde{\Phi}}(\Omega)\}.$$\end{defi}

\nd For the next result we infer the reader to \cite[Lemma 4.1]{Tienari}.

\begin{lem} For any $\epsilon\in (0,1]$ we have
\label{domain}\begin{itemize}
\item[\rm{(i)}] $(1-\epsilon)\mathcal{L}_\Phi(\Omega)\subset \operatorname{dom}(\phi(t)t)$,
\item[\rm{(ii)}] $(1-\epsilon)\mathcal{L}_\Phi(\Omega)+E_\Phi(\Omega)\subset (1-\epsilon/2)
\mathcal{L}_\Phi(\Omega)\subset \operatorname{dom}(\phi(t)t)$.
\end{itemize}
\end{lem}

\nd \proof. At forst we show  $(i)$:  take $u\in \mathcal{L}(\Omega)$ and $\epsilon\in(0,1)$. By the Young inequality (\ref{young}), we easily see that (recall that by taking $s=t\phi(t)$ we obtain the equality in the Young inequality) $\tilde{\Phi}(t\phi(t))\le t(t\phi(t))$. We also see that (due the fact that $t\phi(t)$ is increasing) that $t(t\phi(t))\le \Phi(2t)$.
\vskip.2cm
\nd Therefore \begin{equation}\label{d1}\frac{\epsilon}{1-\epsilon}\tilde{\Phi}((1-\epsilon)u\phi((1-\epsilon)u))\le \epsilon (1-\epsilon)u\phi((1-\epsilon)u).
\end{equation}

\nd Now due the fact that $t\phi(t)$ is increasing, we obtain that
\begin{equation}\label{d2}\epsilon (1-\epsilon)u\phi((1-\epsilon)u)\le \int_{(1-\epsilon)u}^us\phi(s)ds\le \Phi(u).
\end{equation}

\nd Using (\ref{d1}) and (\ref{d2}) together, we conclude the proof of item $i)$.

\nd For the item $ii)$ we put $v\in E_{\Phi}(\Omega)$. Note that $$\frac{v}{1-\epsilon/2}=\frac{2}{\epsilon}\left(1-\frac{1-\epsilon}{1-\epsilon/2}\right)v.$$

\nd Therefore, using the convexity of $\Phi$, we obtain that $$\Phi\left(\frac{1}{1-\epsilon/2}((1-\epsilon)u+v)\right)\le \frac{1-\epsilon}{1-\epsilon/2}\Phi(u)+\left(1-\frac{1-\epsilon}{1-\epsilon/2}\right)\Phi\left(\frac{2v}{\epsilon}\right).$$
\nd This assertion implies that $$(1-\epsilon)\mathcal{L}_\Phi(\Omega)+E_\Phi(\Omega)\subset(1-\epsilon/2)\mathcal{L}_\Phi(\Omega).$$

\nd Now using item $i)$ the proof for item $ii)$ is now finished. \hfill \fbox \hsf

\vskip.2cm
\nd Now we shall prove, with some adapted ideas from  Tienari \cite{Tienari}, that $|\nabla v_k|\in \operatorname{dom}(\phi(t)t)$.

\begin{lem}\label{nablav} The function $|\nabla v_k|$ satisfies $\phi(|\nabla v_k|)|\nabla v_k|\in \mathcal{L}_{\widetilde{\Phi}}(\Omega)$. \end{lem}

\nd \proof. Consider the following sequence $f_k(\epsilon)=I_k(\lambda, (1-\epsilon)v_k)$ for any $\epsilon \in[0,1]$. Recall that $v_k$ is the minimizer of $I_k(\lambda,\cdot)$. This ensures that $f_k(\epsilon)<\infty$ is finite for any $\epsilon\in [0,1]$.  Moreover, we observe that $|\nabla v_k|\in \mathcal{L}_\Phi(\Omega)$. Therefore, using the item i) of Lemma \ref{domain}, we already conclude that $f_k$ is differentiable in the interval $(0,1]$. Furthermore, we know that
\begin{equation}\label{funbounded}
f_k'(\epsilon)=-\int_\Omega \phi((1-\epsilon)|\nabla v_k|)(1-\epsilon)|\nabla v_k|^2 dx +\lambda\int_\Omega f_k((1-\epsilon)v_k)v_k dx.
\end{equation}

\nd Suppose, on the contrary that $\phi(|\nabla v_k|)|\nabla v_k|\notin \mathcal{L}_{\tilde{\Phi}}(\Omega)$. Recall from Young inequality ($\ref{young}$) that
$$
\int_\Omega \phi(|\nabla v_k|)|\nabla v_k|^2 dx =\int_\Omega \Phi(|\nabla v_k|)dx +\int_\Omega \tilde{\Phi}(\phi|\nabla v_k|)|\nabla v_k| dx =\infty.
$$

\nd From the last equality and ($\ref{funbounded}$), we mention that $f'(\epsilon)\to -\infty$ if $\epsilon\to 0$. Hence, there exists $\epsilon_0\in (0,1]$ such that $f(\epsilon_0)<f(0)$ which is an contradiction because of $v_r$ is a minimizer for $I_k(\lambda,\cdot)$.
\hfill \fbox \hsf
\vskip.3cm

\nd With the help of the last Lemma, we can now prove that the minimum we found in Lemma \ref{Qk not-empty} does satisfies the Euler equation:
\begin{lem}\label{weaksol} For $v_k$ as defined in lemma \eqref{Qk not-empty} we have
	$$
	\int_\Omega \phi(|\nabla v_k|)\nabla v_k\nabla w dx =\lambda \int_\Omega f(v_k)w dx,\  w\in W_0^1 L_{\Phi}(\Omega).
	$$
	\end{lem}

\nd \proof. Let $v\in W_0^1 E_\Phi(\Omega)$ and define $$f_v(\epsilon)=I_k(\lambda,v_\epsilon),\  \epsilon \in [0,1],$$

\nd where $v_\epsilon=1/(1-\epsilon/2)[(1-\epsilon)v_k+\epsilon v]$. The item (ii) of lemma \ref{domain} and Young inequality ($\ref{young}$) imply that $f_v(\epsilon)<\infty $ for all $\epsilon\in [0,1]$. Once $v_k$ is a minimizer of $I_k(\lambda,\cdot)$, we have that \begin{equation}\label{blah}
0\le \frac{f_v(\epsilon)-f_v(0)}{\epsilon},\  v\in W_0^1 E_\Phi(\Omega).
\end{equation}

\nd Note that for $\epsilon<2/3$, the monotonicity of $\phi(t)t$ and the triangle inequality, the following is true
\begin{eqnarray*}
          \left|\frac{\Phi(|\nabla v_\epsilon|)-\Phi(|\nabla v_k|)}{\epsilon}\right| &\le& (\phi(|\nabla v_\epsilon|)|\nabla v_\epsilon|+\phi(|\nabla u_k|)|\nabla u_k|)\frac{|\nabla v_\epsilon-\nabla v_k|}{\epsilon}, \\
             &\le& (2\phi(|\nabla v_k|)|\nabla v_k|+\phi(|\nabla v|)|\nabla v|)(|\nabla v_k|+|\nabla v|).
          \end{eqnarray*}

\nd As $|\nabla v|\in \operatorname{dom}(\phi(t)t)$ for every $v\in W^1_0 E_\Phi(\Omega)$ (this is true because of the inequality $\tilde{\Phi}(\phi(t)t)\le \Phi(2t)$), we conclude that the right hand side of the last inequality, is a function in $L^1(\Omega)$. Now using the fact that $$\frac{\Phi(|\nabla v_\epsilon|)-\Phi(|\nabla v_k|)}{\epsilon}\to \phi(|\nabla v_k|)\nabla v_k(\nabla v-\nabla v_k/2),\ \mbox{a.e.},\ x\in \Omega,\ \mbox{when}\ \epsilon\to 0,$$

\nd we infer from inequality ($\ref{blah}$) and Lebesgue theorem that
$$0
\le \int_\Omega \phi(|\nabla v_k|)\nabla v_k(\nabla v-\nabla v_k/2) dx-\lambda \int_\Omega f(v_k)(v-v_k/2) dx,\  v\in W_0^1 E_\Phi(\Omega).
$$

\nd As a consequence
$$
\int_\Omega \phi(|\nabla v_k|)\nabla v_k\nabla v=\lambda \int_\Omega f(v_k)v dx,\ l v\in W_0^1 E_\Phi(\Omega).
$$

\nd At this moment, using the weak star density of $W_0^1 E_\Phi(\Omega)$ in $W_0^1 L_\Phi(\Omega)$, we mention that

$$
\int_\Omega \phi(|\nabla v_k|)\nabla v_k\nabla v dx=\lambda \int_\Omega f(v_k)v dx,\   v\in W_0^1 L_\Phi(\Omega).
$$
Hence $v_k$ is a weak solution to the problem \eqref{3.1}. This finishes the proof for this lemma. \hfill \fbox \hsf

 \nd To continue, we shall prove that $I_{k}(\lambda, \cdot)$ admits at least one weak solution $v_{k}$ satisfying $a_{k} < \|v_{k}\|_{\infty} \leq a_{k+1}, k = 1, \ldots, m-1$. This is crucial in order to get our main result.

\begin{lem}\label {lema 2.3}
  There is $\lambda_k>0$ such that
$$
a_{k-1}<\|v_{k}\|_\infty \leq a_{k}
$$
\nd for each minimum $v_k \equiv v_k(\lambda)$  of $I_{k}(\lambda, \cdot) $ with $\lambda > \lambda_k$.
\end{lem}

\nd \proof. The proof is similar to those given in \cite{Hess} and \cite{loc-schmitt} and the idea is the following: Take $\delta > 0$ and consider the open set
 $$
\Omega_\delta=\{x\in\Omega~|~ \operatorname{dist}(x,\partial\Omega)<\delta\}.
$$
\nd Set
$$
{\widetilde{\alpha }}_k :=F(a_k)-\max{\{F(s)~|~ 0\leq s\leq a_{k-1}\}}
$$
\nd and note that by $(f_3)$, $ {\widetilde{\alpha }}_k > 0$. Choose  $w_\delta\in C_0^\infty(\Omega)$ such that
$$
0\leq w_\delta\leq a_k~\mbox{and}~ w_\delta(x)=a_k,~  x\in\Omega\setminus\Omega_\delta.
$$
\nd Writing $\Omega=\Omega_{\delta} \cup (\Omega \setminus \Omega_{\delta})$ and setting $ C_k=\max{\{|F(s)|~|~\ 0\leq s\leq a_k\}}$ we get to,

$$
 \int_\Omega F(w_\delta) dx \geq \int_\Omega F(a_k) dx -2C_k|\Omega_\delta|. \nonumber
$$

\nd Let  $u\in W_0^1L_\Phi(\Omega)$ such that $0\leq u\leq a_{k-1}$. By the  inequality above we have
$$
\int_\Omega F(w_\delta) dx -\int_\Omega F(u) dx \geq {\widetilde{\alpha }}_k |\Omega|-2C_k|\Omega_\delta|.
$$
\nd Since  $|\Omega_\delta|\to 0$ as $\delta\to 0$ there is $\delta>0$ such that
$$  \eta_k := {\widetilde{\alpha }}_k |\Omega|-2C_k|\Omega_\delta|>0.
$$
\nd Set $w=w_\delta$ and pick $u\in W_0^1L_\Phi(\Omega)$ with $0\leq u\leq a_{k-1}$. Choosing $\lambda_k > 0$ large enough, taking $\lambda \geq \lambda_k$ and  making use of the expessions of  $I_k(\lambda,w_\delta), I_{k-1}(\lambda,u)$ and the inequality just above we infer that
\begin{equation}
 I_k(\lambda,w_\delta)-I_{k-1}(\lambda,u)  \leq \int_\Omega\Phi(|\nabla w_\delta|) dx -\lambda\eta_k
   < 0
\end{equation}
\nd and  hence
\begin{equation}\label{1.4}I
_k(\lambda,w_\delta)<I_{k-1}(\lambda,u)~\mbox{ for }~\lambda \geq \lambda_k.
\end{equation}

\nd To finish, assume, on  the contrary,  that there is a minimum $v_k(\lambda)$ of  $I_k(\lambda,\cdot)$  such that $v_k(\lambda) \leq a_{k-1}$. It follows by $(\ref{1.4})$ and lemma \ref{lema 2.1} that
$$
I_k(\lambda,w_\delta)<I_{k-1}(\lambda,v_k(\lambda)).
$$
\nd   On the other hand, because $I_{k-1}(\lambda,v_k(\lambda))=I_k(\lambda,v_k(\lambda))$ and since
$v_k(\lambda)$ is a minimum of $I_k(\lambda,\cdot)$ we have
$$
 I_{k-1}(\lambda,v_k(\lambda))=I_k(\lambda,v_k)  \leq  I_k(\lambda,w_\delta) .     \nonumber \\
  $$
\nd  The inequalities just above  lead to a contradiction.\hfill \fbox \hsf


\section{Proof of Theorem \ref{1.1} } The proof is based on Loc \& Schmitt \cite{loc-schmitt}.  However, we will get into  details   taking into account the Orlicz-Sobolev  spaces framework. In this sense we will make use of a general result on lower and upper solutions by Le \cite[theorem 3.2]{Khoi}.
\vskip.1cm

\nd For the proof of {\rm(ii)} we will need the lemma below. In order to state it, take an open ball $B$ centered at $0$ with radius $R$ containing $\Omega$. Consider the functions  $\alpha, \beta : {\overline{B}} \to \r$ defined as follows:
$$
 \alpha(x)= \left\{ \begin{array}{rll}
 u(x), &~ x \in  \overline{\Omega}  \\
  0, &~  x \in \overline{B} \setminus \Omega,
       \end{array} \right.
~~~~~ \beta(x)=a_{k+1},~   x \in \overline{B}.
$$
\nd Since $u \in  W_0^1L_\Phi(\Omega)$, we know by the Lemma (\ref{LA1}) in the Appendix, that the extension by zero $\overline{u}$ of $u$ belongs to $W_0^1L_\Phi(\mathbb{R}^N)$. Hence, due to the fact that the extension by zero of the function $\alpha$ to $\mathbb{R}^N$, coincides with $\overline{u}$, we also conclude by using Lemma (\ref{LA1}) (appendix) again that $\alpha\in W_0^1L_\Phi(B)$.

\begin{lem}\label{lemma 3.2 Kaye} The functions
$\beta$ and $\alpha$ are respectively upper and lower solutions to the elliptic problem
\begin{equation}\label{3.2}
 \left\{ \begin{array}{cl}
 -\Delta_\Phi u=\lambda f(u)~\mbox{in}~ B, \\
  \\
    u\in W_0^1L_\Phi(B).
       \end{array} \right.
\end{equation}
\end{lem}

\nd \proof. That $\beta$ is an upper solution is immediately. Now we shall prove that $\alpha$ is a subsolution for the problem \eqref{3.2}. Let $v_n(x)=n\min\{u(x),1/n\}$ be a fixed function. Note that $u_n(x)\to 1$ for each $x$ where $u(x)\neq 0$ and $u(x)\to 0$ for each $x$ where $u(x)=0$. Moreover, on the set $\{x\in \Omega:\ u(x)=0\}$, we have that $\nabla u(x)=0$ a.e..  Therefore, for each $v\in W_0^1L_\Phi(\Omega)$ with $v\ge 0$, the following inequalities are true  \begin{eqnarray*}
                            \int_B \phi(|\nabla \alpha|)\nabla \alpha\nabla vdx  &=& \int_\Omega \phi(|\nabla u|)\nabla u\nabla v dx \\
                             &=& \lim_{n\to \infty} \int_\Omega v_n\phi(|\nabla u|)\nabla u\nabla v dx \\
                             &=& \lim_{n\to\infty} \left(\int_\Omega \phi(|\nabla u|)\nabla u\nabla (v_nv) dx-\int_\Omega v\phi(|\nabla u|)\nabla u\nabla v_ndx\right)  \\
                             &\le& \lim_{n\to \infty}\int_\Omega f(u)v_nv dx \\
                             &\le & \int_B f(\alpha)vdx.
                          \end{eqnarray*}
This proves the lemma.
\hfill \fbox \hsf

\nd {\bf Proof of  {\rm(i)} of theorem \ref{1.1}}. Take  $\lambda > 0$. By lemma  \ref{Qk not-empty},   for each  $k = 2, \cdots, m$ there is a minimum $v_k \equiv v_k(\lambda)$ of $I_k(\lambda)$,  which is actually a weak solution of problem (\ref{3.1k}).  By lemma \ref{lema 2.1}, we know that $0 \leq v_k \leq a_{k}~\mbox{a.e. in}~\Omega$.
\vskip.2cm

\nd Now we mention that, using Lemma \ref {lema 2.3}, there is $\overline{\lambda} \geq \displaystyle \max_{2 \leq k \leq m} \{\lambda_k \}$ such that $v_2, \cdots, v_m$ are  solutions of problem (\ref{3.1}) for any $\lambda > \overline{\lambda}$. Moreover, these solutions satisfy
$$
a_1 <  \|v_{2}\|_\infty \leq a_2 < \|v_{3}\|_\infty \leq \cdots \leq a_{m-1}<\|v_{m}\|_\infty \leq a_{m}
$$
\nd Now we take $u_{k-1}  \equiv v_{k}(\lambda),~k = 2, \cdots, m$. This ends the proof of the first part of theorem \ref{1.1} .
\vskip.2cm

\nd {\bf Proof of  {\rm(ii)} of theorem \ref{1.1}}. We distinguish between two cases.
\vskip.3cm

\nd {\bf Case 1}~~ $f(0) > 0$.
\vskip.2cm

\nd This case is  more difficult. In order to address it we state and prove the lemma below.
\begin{lem}\label{thm 3.1 Kaye}
Assume $(\phi_1)-(\phi_4)$, $(f_1)-(f_2)$ and $f(0)>0$. If $u$ is a non-negative weak solution of $(\ref{1.1})$ such that  $a_{k-1} < \|u\|_\infty \leq a_{k}$ then
 $$\int_{a_k}^{a_{k+1}}f(s) ds>0.$$
\end{lem}

\nd \proof~  Let us consider the case that $k = 2$, the other cases may be treated in a similar manner. Take the lower and upper solutions respectively $\alpha$ and $a_2$ for (\ref{3.2}).
\vskip.1cm

\nd Now, applying Theorems 3.2, 4.1 and 5.1 of \cite{Khoi}, we find a maximal solution say $\overline{u}$ of (\ref{3.2}) such that $\alpha (x)\leq\overline{u}(x)\leq a_{2}$ for $x\in B$.
\vskip.2cm

\nd The verification of  the {\bf Claim} below follows as in \cite{loc-schmitt}.
 \begin{claim}\label{rad symm}
 {\it $\overline{u}$ is radially symmetric, i.e.  $\overline{u}(x_1)=\overline{u}(x_2),~x_i \in B,~ |x_1|=|x_2|$}.
\end{claim}

\vskip.1cm

\nd Now we set
$$
u(r) = \overline{u}(x)~\mbox{where}~r = |x|~\mbox{and}~  x \in B,
$$
\nd and because the extension by zero outside of $\Omega$, of the function $\overline{u}$, is an absolutely continuous function, with respect to a.e. segment of line in the direction of a vector $\eta\in \{y\in \mathbb{R}^N:\ \|y\|=1\}$, we conclude that $u$ is continuous in $(0,R]$ and
$$
u\in W^{1,1}(0,R),\ u(1)=0.
$$
\nd  Let $r \in (0,R)$ and pick $\epsilon > 0$ small such that  $r + \epsilon < R$. Note that
\begin{equation}\label{Eq u bar}
\int_B \phi(|\nabla\overline{u}|) \nabla \overline {u}\nabla v  dx =\lambda\int_B f(\overline{u}) v dx,~  v\in W_0^{1, \Phi}({B}).
\end{equation}
\nd  Consider the radially symmetric cut-off function $v_{r,\epsilon}(x) = v_{r,\epsilon}(r) $, where
$$
v_{r,\epsilon}(t) := \left\{ \begin{array}{l}
1~~ \mbox{if}~~ 0 \leq t \leq r,\\
linear~~ \mbox{if}~~ r \leq t \leq r + \epsilon,\\
0~~ \mbox{if}~~ r+\epsilon \leq t \leq R.
\end{array} \right.
$$
\nd and notice that  $v_{r,\epsilon} \in W_0^{1, \Phi}({B}) \cap Lip( {\overline{B}}    )$. Setting  $v = v_{r,\epsilon}$  in (\ref{Eq u bar})  and using the radial symmetry we get to
$$
\frac{-1}{\epsilon} \int_{r}^{r + \epsilon} t^{N-1}
\phi(| {{u}}^{\prime}|)
 {{u}}^{\prime}~dt =
\int_0^r t^{N-1}~ \lambda
f({u})~dt + \int_{r}^{r+\epsilon} t^{N-1} \lambda f({u})\upsilon dt. \\
~~
$$
\nd Once $t^{N-1}\phi(|u'|)u'\in L^1(0,R)$, we conlude by letting $\epsilon \to 0$, that for a.e. $r\in (0,R)$
\begin{equation}\label{Eq u bar limit}
- r^{N-1} \phi(| {{u}}^{\prime}(r)|) {{u}}^{\prime}(r) =
\displaystyle \int_0^r ~ \lambda
f({u}) t^{N-1}~dt,~0 < r < R.
\end{equation}
\nd From equation (\ref{Eq u bar limit}) and by the boundedness of $u$, we obtain that $u'$ is continuous and $u'(0)=0$. Set
$$
\|u\|_\infty = \max \{ u(r)~|~r \in [0,R] \},
$$
\nd and choose numbers $r_0, r_1 \in [0,R)$ with $r_1 \in (r_0, R)$ such that
$$
u(r_0) = \|u \|_{\infty}~\mbox{and}~ u(r_1) = a_1.
$$
\nd Note that
$$
u(r_0)  >  u(r_1)~\mbox{and}~ 0 \leq r_0 < r_1 < R.
$$
\begin{claim}\label{b1}
~~~~  $\| u \|_{\infty} > b_1$.
\end{claim}
\nd Indeed, assume  on the contrary that,  $u(r_0) \leq b_1$. Take $\delta > 0$ small such that
$$
a_1 < u(r) \leq u(r_0),~~ r_0  \leq  r \leq   r_0 + \delta.
$$
\nd We have by (\ref{Eq u bar limit})
\begin{equation}\label {r0}
-r_{0}^{N-1} \phi(| {{u}}^{\prime}(r_0)|) {{u}}^{\prime}(r_0) = \displaystyle \int_{0}^{r_0} ~ \lambda
f({u}) t^{N-1}~dt,
\end{equation}
\begin{equation}\label{rr}
-r^{N-1} \phi(| {{u}}^{\prime}(r)|) {{u}}^{\prime}(r) = \displaystyle\int_{0}^{r} ~ \lambda
f({u}) t^{N-1}~dt.
\end{equation}
\nd Subtracting (\ref{rr}) from (\ref{r0}) and recalling that  $u^{\prime}(r_0) = 0$ one obtains,
$$
-r^{N-1} \phi(| {{u}}^{\prime}(r)|) {{u}}^{\prime}(r) = \int_{r_0}^{{r}} ~ \lambda
f({u}) t^{N-1}~dt,~  r_0  \leq  r \leq   r_0 + \delta.
$$
\nd Since $f \leq 0~\mbox{on}~[a_1,b_1]$,
$$
r^{N-1} \phi(| {{u}}^{\prime}(r)|) {{u}}^{\prime}(r) \geq 0,~ r_0 \leq r \leq r_0 + \delta.
$$
\nd It follows that ${{u}}^{\prime}(r) \geq 0~\mbox{for}~ r_0 \leq r \leq r_0 + \delta$. But, since $u(r_0)$ is a global maximum on $[0,R]$, it follows that $u^{\prime} = 0~\mbox{on}~ [ r_0,  r_0 + \delta]$. By a continuation argument we get $u^{\prime} = 0~\mbox{on}~ [r_0, r_1)$ so that $u = \| u \|_{\infty}~ \mbox{on}~ [r_0, r_1]$, contradicting $u(r_0) > a_1$. As a consequence,  $\| u \|_{\infty}  >  b_1$, proving  Claim \ref{b1}.

\begin{claim}
~~ $u \in C^{2}({\cal{O}})$~where ${\cal{O}} := \{r \in (0,R)~|~ u'(r) \neq 0 \}$.
\end{claim}

\nd Of course ${\cal{O}}$  is an open set. Motivated by the left hand side of (\ref{Eq u bar limit}) consider
$$
G(z)  =   \phi(z) z,~~ z \in \r ,
$$
\nd where $z$ is set to play the role of $u^{\prime}$.  Recall that
$$
G~\mbox{is odd},~~  G^{\prime}(z) = (\phi(z) z )^{\prime} > 0~\mbox{for}~z > 0
$$
\nd and
$$
G(z) =  \phi(| {{u}}^{\prime}(r)|) {{u}}^{\prime}(r) =  -\frac{1}{r^{N-1}}\int_{0}^{{r}} ~ \lambda
f({u}) t^{N-1}~dt.
$$
\nd Since $\phi(z)z \in C^{1}$ and $(\phi(z)z)'\neq 0$ for $z\neq 0$, we get by applying the Inverse Function Theorem  in ${\cal{O}}$ that $z = z(r,u)$ is a $C^{1}$-function of $r$. Since $z = u^{\prime}$, the claim is proved.

\begin{claim}\label{int positive}
   ~~ $\int_{ a_{1}}^{\|u \|_{\infty}} f(s) ds > 0$.
\end{claim}

\nd Differentiating in (\ref{Eq u bar limit}) and multiplying by $u^{\prime}$ we get
$$
\big ( t^{N-1} \phi(| {{u}}^{\prime}(t)|) {{u}}^{\prime}(t)  \big)^{\prime} u^{\prime}(t) = -\lambda
f({u(t)}) u^{\prime}(t) t^{N-1},
$$
\nd and hence
$$
\frac{(N-1)}{t}   \phi(| u^{\prime}|) (u^{\prime})^2 + ( \phi(| {u}^{\prime}|) u^{\prime})^{\prime}  u^{\prime} = -\lambda
f({u}) u^{\prime}.
$$
\nd Integrating from $r_0$ to $r_1$ we have
\begin{equation}\label{before}
-\Big[\int_{r_0}^{r_1}  \frac{(N-1)}{t}   \phi(| u^{\prime}|) (u^{\prime})^2  dt +\int_{r_0}^{r_1} \big[\phi(| {u}^{\prime}| )u^{\prime} \big]^{\prime}  u^{\prime}  dt \Big]  =   \int _{r_0}^{r_1} \lambda f({u}) u^{\prime} dt.
\end{equation}
\nd Making the change of variables  $s = u^{\prime}(t)$ in the second and third integrals in (\ref{before})  and applying the arguments in \cite{loc-schmitt} leads to  {\bf Claim  \ref{int positive}}.
\vskip.2cm

\nd Since $f \geq 0$ on $(b_1,a_2)$ and $\|u\|_{\infty} > b_1$ it follows that
$$
\int ^{a_2}_{a_1}
f({s}) ds  > 0,
$$
\nd ending the proof of  lemma \ref{thm 3.1 Kaye}.  \hfill \fbox \hsf
\vskip.2cm

\nd {\bf Case 2}~~$f(0) = 0$.
\vskip.1cm

\nd This case is handled as in \cite{loc-schmitt}. The Theorem 1.1 is now proved.  \hfill \fbox \hsf
\vskip.2cm

\section {Appendix }

\noindent In this Section we state and prove  a version of the Stampacchia Theorem (generalized Chain Rule) for the case of nonreflexive Orlicz-Sobolev spaces.
\vskip.3cm

\begin{prop}\label{CA1} Let $g:\mathbb{R}\to \mathbb{R}$ be a Lipschitz continuous function such that  $\|g'\|_\infty<M$ and $g(0)=0$. Let $u\in W_0^{1}L_{\Phi}(\Omega)$.   Then  $g(u)\in W_0^{1}L_{\Phi}(\Omega)$ and
$$
 \frac{\partial g(u)}{\partial x_i} = g^{\prime}(u)  \frac{\partial u}{\partial x_i}~~\mbox{a.e. in}~~ \Omega.
$$
\end{prop}

\nd At first we recall, for the reader's convenience, the definition and basic properties of the trace on $\partial \Omega$ of an element of $W^1L_\Phi(\Omega)$.
 We refer the reader to Gossez \cite{Gz3}.
\vskip.2cm

\noindent Let $\gamma:C^\infty(\overline{\Omega})\to C(\partial \Omega)$ be defined by the linear map $\gamma(u) = u_{| \partial \Omega}$. Then  $\gamma$ is  continuous with respect to the topologies
$$
\sigma\left(\prod L_\Phi(\Omega),\prod E_{\widetilde{\Phi}}(\Omega)\right)~ \mbox{and}~ \sigma (L_\Phi(\partial \Omega),E_{\tilde{\Phi}}(\partial\Omega)).
$$

\nd Using the facts that  $C^\infty(\overline{\Omega})$ is dense in $W^1L_\Phi(\Omega)$ with respect to the topology $\sigma\left(\prod L_\Phi(\Omega),\prod E_{\widetilde{\Phi}}(\Omega)\right)$ and $C(\partial \Omega)$ is dense in $L_{\Phi}(\partial \Omega)$ with respect to the topology $\sigma (L_\Phi(\partial \Omega),E_{\tilde{\Phi}}(\partial\Omega))$, $\gamma$ admits an only continuous extension to a linear map namely  $\gamma:W_0^1L_\Phi(\Omega)\to L_{\Phi}(\partial\Omega)$.

\nd It can be shown that
$$
W_0^1L_\Phi(\Omega)=\{u\in W^1L_\Phi(\Omega)~|~ \gamma(u)=0\}.
$$

\nd Let $u:\Omega\to\r$. We define
$\overline{u}:\mathbb{R}^N\to \mathbb{R}$ by
$$
 \overline{u}(x)= \left\{ \begin{array}{rll}
 u(x) &\mbox{if}~ x \in  {\Omega},  \\
  0 & \mbox{if}~  x \in  {\Omega}^c.
       \end{array} \right.
$$

\noindent Let $\nu=(\nu_1,\cdots,\nu_N)$ be the outward unit normal vector field of $\partial \Omega$. The Green's formula reads as,
 \begin{equation}\label{GI}\int_\Omega u\frac{\partial v}{\partial x_i}dx +\int_\Omega v\frac{\partial u}{\partial x_i}dx =\int_{\partial \Omega}\gamma (u)\gamma (v)\nu_i dx,~i=1,\cdots, N,
 \end{equation}
\nd where   $u\in W^1L_\Phi(\Omega)$, $v\in W^1L_{\widetilde{\Phi}}(\Omega)$.
\vskip.2cm

\nd We will give the proof of the result below:

\begin{lem}\label{LA1}
	$$
	W_0^{1}L_{\Phi}(\Omega)=\{u\in W^{1}L_{\Phi}(\Omega)~|~ \overline{u}\in W^{1}L_{\Phi}(\mathbb{R}^N)\}.
	$$
	\end{lem}
	
\nd \proof~ Set
$$
E = \{u\in W^{1}L_{\Phi}(\Omega)~|~ \overline{u}\in W^{1}L_{\Phi}(\mathbb{R}^N)\}.
$$	
Take $u\in W_0^{1}L_{\Phi}(\Omega)$ and let $g_i:{\r}^N\to \mathbb{R}$ be defined by
$$
 g_i(x) = \left\{ \begin{array}{rll}
 (\partial u/\partial x_i) (x) &\mbox{if}~ x \in  {\Omega},  \\
  0 & \mbox{if}~  x \in  {\Omega}^c.
       \end{array} \right.
$$
\nd Computing derivatives in the distribution sense and using the generalized Green's formula we have
\begin{eqnarray*}
         \int_{{\r}^N} u\frac{\partial\varphi}{\partial x_i}dx &=& \int_\Omega u\frac{\partial\varphi}{\partial x_i}dx, \\
             &=& -\int_\Omega\frac{\partial u}{\partial x_i}\varphi dx, \\
             &=&  -\int_{{\r}^N} g_i\varphi dx,\ \varphi\in C_0^\infty ({\r}^N).
          \end{eqnarray*}
\nd Using the fact that $\overline{u},g_i\in L_{\Phi}(\mathbb{R}^N)$, we conclude that $u\in E$.  On the other hand, take $u\in E$. We claim that $\gamma(u)=0$. Indeed, by the definition of weak derivative and Green's identity it follows  that for all $v\in W^{1,\tilde\Phi}(\Omega)$ and $i=1,\cdots,N$, we have
$$
\int_{\partial \Omega} \gamma(u)\gamma(v)\nu_i dx=0.
$$
\nd and so
$$
\int_{\partial \Omega} \gamma(u)\nu_i w dx=0, \  w\in C(\partial \Omega).
$$
\nd Because $\gamma (u)\nu_i\in L_{\Phi}(\partial \Omega)$ and $C(\partial\Omega)$ is  dense (with respect to the norm topology) in $E_{\widetilde{\Phi}}(\partial \Omega)$, we conclude that
$$
\int_{\partial \Omega} \gamma(u)\nu_i w dx = 0, \  w\in E_{\widetilde{\Phi}}(\partial \Omega).
$$

\nd By taking $w(x)=1$ if $\gamma(u)\nu_i$ is positive, $w(x)=0$ if $\gamma(u)\nu_i$ is zero and $w(x)=-1$ if $\gamma(u)\nu_i$ is negative, we conclude that $w\in L^\infty(\partial\Omega)\subset E^{\tilde\Phi}(\partial \Omega)$ and
$$
\int_{\partial \Omega} |\gamma(u)\nu_i| d \sigma=0.
$$
\nd So, $\gamma (u)=0$ a.e. on the supoort of $\nu_i$ and because $$\bigcup_{i=1}^N \operatorname{supp}(\nu_i)=\partial \Omega,$$

\nd we conclude that $\gamma(u)=0$ and hence $u\in W_0^{1}L_{\Phi}( \Omega)$. This finishes the proof. \hfill \fbox \hsf
\vskip.3cm

\nd \proof~ {\bf of Proposition \ref{CA1}.} Indeed, by lemma \ref{LA1},  $\overline{u}\in W^{1}L_{\Phi}(\mathbb{R}^N)$.
Reminding that
\begin{equation}\label{emb-1}
W^{1}L_{\Phi}(\Omega) \hookrightarrow W^{1,1}(\Omega),
\end{equation}
\nd we have  $\overline{u}\in W^{1,1}(\mathbb{R}^N)$ and by the Chain rule for  $W^{1,1}$, (cf. Gilbarg \& Trudinger \cite{GT}), we infer that $g(\overline{u})\in W^{1,1}(\mathbb{R}^N)$ and in addition,
$$
\int_{\Omega}  g(u) \frac{\partial\varphi}{\partial x_i}dx = \int_{\Omega } g^{\prime}(u)  \frac{\partial u}{\partial x_i} \varphi dx,\ \varphi\in C_0^\infty (\Omega).
$$
\nd We will show next that $\overline{g(u)} \in W^{1}L_{\Phi}(\mathbb{R}^N)$. Indeed, on one hand, we have that $g(\overline{u})=\overline{g(u)}$. On the other hand, setting
 $$
 h_i(x) = \left\{ \begin{array}{rll}
 \frac{\partial g(u)}{\partial x_i}(x) &\mbox{if}~ x \in  {\Omega},  \\
  0 & \mbox{if}~  x \in  {\Omega}^c,
       \end{array} \right.
 $$
\nd  we have $\overline{g(u)},h_i\in L^{\Phi}(\mathbb{R}^N)$. Therefore, $\overline{g(u)}\in W^{1}L_{\Phi}(\mathbb{R}^N)$. Now, using lemma \ref{LA1} again, we obtain that $g(u)\in W_0^{1}L_{\Phi}(\Omega)$  and
$$
 \frac{\partial g(u)}{\partial x_i} = g^{\prime}(u)  \frac{\partial u}{\partial x_i}~\mbox{a.e. in}~ \Omega.
$$
\nd This concludes the proof of Proposition \ref{CA1}. \hfill \fbox \hsf

\vskip.2cm

\begin{flushright}
\scriptsize{ Edcarlos D. Silva \\

 Jos\'e V.A. Goncalves\\

Kaye O. Silva}\\
\medskip

  \scriptsize{Universidade Federal de Goi\'as\\
   Instituto de Matem\'atica e Estat\'istica\\
   74001-970 Goi\^ania, GO - Brasil}\\
\vskip.2cm

\scriptsize{ Emails: edcarlos@ufg.br\\

goncalves.jva@gmail.com\\

kayeoliveira@hotmail.com\\

}
\end{flushright}

\end{document}